\documentclass[12pt,a4paper]{amsart}
\usepackage{amscd}
\usepackage{mathrsfs}
\usepackage{amsfonts}
\usepackage{graphicx}
\usepackage{amsmath,amscd,stmaryrd,amssymb,array, amsfonts,amsmath,amssymb}
\usepackage{enumitem}
\usepackage{amsthm}
\usepackage{bm}

\usepackage{hyperref}

\setenumerate[1]{itemsep=0pt,partopsep=0pt,parsep=\parskip,topsep=5pt}
\setitemize[1]{itemsep=0pt,partopsep=0pt,parsep=\parskip,topsep=0pt}
\setdescription{itemsep=0pt,partopsep=0pt,parsep=\parskip,topsep=5pt}

\numberwithin{figure}{section}
 \numberwithin{equation}{section}

\newtheorem{theorem}{Theorem}[section]
\newtheorem{proposition}[theorem]{Proposition}
\newtheorem{definition}[theorem]{Definition}

\newtheorem{lemma}[theorem]{Lemma}
\newtheorem{remark}[theorem]{Remark}

\newtheorem{claim}[theorem]{Claim}

\def\be{\begin{equation}}
\def\ee{\end{equation}}
\def\bes{\begin{equation*}}
\def\ees{\end{equation*}}
\def\bsp{\begin{split}}
\def\esp{\end{split}}

\def\ba{\begin{array}}
\def\ea{\end{array}}
\def\benu{\begin{enumerate}}
\def\eenu{\end{enumerate}}
\def\bt{\begin{theorem}}
\def\et{\end{theorem}}
\def\bp{\begin{proposition}}
\def\ep{\end{proposition}}
\def\bl{\begin{lemma}}
\def\el{\end{lemma}}
\def\br{\begin{remark}}
\def\er{\end{remark}}
\def\bd{\begin{definition}}
\def\ed{\end{definition}}

\def\.{\cdot}

\def\~{\tilde}
\def\8{\infty}

\def\Vs{\vskip8pt}\def\vs{\vskip4pt}

\def\({\left(}\def\){\right)}

\begin{document}

\begin{center}
{\bf\Large
A Semilinear Elliptic Problem with Critical Exponent and Potential Terms}
\end{center}

\vs\centerline{Haoyu Li}  
\begin{center}
{\footnotesize
{Departamento de Matem\'atica,  Universidade Federal de S$\~{a}$o Carlos,\\
          S$\~{a}$o Carlos-SP, 13565-905, Brazil\\
{\em E-mail}:  e-mail: hyli1994@hotmail.com}}
\end{center}
\vs\centerline{Li Ma}
\begin{center}
{\footnotesize
{School of Mathematics and Physics, University of Science and Technology Beijing, \\
30 Xueyuan Road, Haidian District, Beijing 100083, China\\
and\\
Department of Mathematics, Henan Normal University,\\
           Xinxiang 453007, China\\
{\em E-mail}:  lma17@ustb.edu.cn}}
\end{center}

\Vs

{\footnotesize
{\bf Abstract.}
This paper addresses the following problem.
\begin{equation}
    \left\{
   \begin{array}{lr}
     -{\Delta}u=\lambda I_\alpha*_\Omega u+|u|^{2^*-2}u\mbox{ in }\Omega ,\nonumber\\
     u\in H_0^1(\Omega).\nonumber
   \end{array}
   \right.
\end{equation}
Here, $\Omega$ is a bounded domain in $\mathbb{R}^N$ with $N\geq3$, $2^*=\frac{2N}{N-2}$, $\lambda\in\mathbb{R}$, $\lambda\in(0,N)$, $I_\alpha$ is the Riesz potential and
\begin{align}
I_\alpha*_\Omega u(x):=\int_\Omega \frac{\Gamma(\frac{N-\alpha}{2})}{\Gamma(\frac{\alpha}{2})\pi^\frac{N}{2}2^\alpha|x-y|^{N-\alpha}} u(y)dy. \nonumber
\end{align}
We study the non-existence, existence and multiplicity results. Our argument combines Brezis-Nirenberg's method with the regularity results involving potential terms. Especially, we study the following nonlocal eigenvalue problem.
\begin{equation}
    \left\{
   \begin{array}{lr}
     -{\Delta}u=\lambda I_\alpha*_\Omega u\mbox{ in }\Omega ,\nonumber\\
     \lambda\in\mathbb{R},\,u\in H_0^1(\Omega).\nonumber
   \end{array}
   \right.
\end{equation}

 \Vs
\Vs {\bf 2020 MSC:} 35J20, 35B33, 35J61, 35J61, 45K05, 47A75, 35B50}

\tableofcontents

\section{Introduction}
\subsection{Introduction and historical remarks}
In this paper, we consider the following problem.
\begin{equation}\label{e:000}
    \left\{
   \begin{array}{lr}
     -{\Delta}u=\lambda I_\alpha*_\Omega u+|u|^{2^*-2}u\mbox{ in }\Omega ,\\
     u\in H_0^1(\Omega).
   \end{array}
   \right.
\end{equation}
Here, $\Omega$ is a bounded domain in $\mathbb{R}^N$ with smooth boundary $\partial\Omega$ and $N\geq3$, $2^*=\frac{2N}{N-2}$, $I_\alpha$ is the Riesz potential defined as
\begin{align}
I_\alpha=\frac{\Gamma(\frac{N-\alpha}{2})}{\Gamma(\frac{\alpha}{2})\pi^\frac{N}{2}2^\alpha|x|^{N-\alpha}} \nonumber
\end{align}
and the corresponding potential term is defined as
\begin{align}
I_\alpha*_\Omega u(x):=\int_\Omega \frac{\Gamma(\frac{N-\alpha}{2})}{\Gamma(\frac{\alpha}{2})\pi^\frac{N}{2}2^\alpha|x-y|^{N-\alpha}} u(y)dy. \nonumber
\end{align}
In follows, we denote $C_{N,\alpha}=\frac{\Gamma(\frac{N-\alpha}{2})}{\Gamma(\frac{\alpha}{2})\pi^\frac{N}{2}2^\alpha}$ for the sake of brevity.
Our topic is stream out of the following problem involving a convolution term.
\begin{equation}\label{e:100}
    \left\{
   \begin{array}{lr}
     -\Delta u+\omega u+\lambda(I_\alpha* F(u))F'(u)=g(x,u)\mbox{ in }\Omega ,\\
     u\in H_0^1(\Omega).
   \end{array}
   \right.
\end{equation}
The equation in the form of (\ref{e:100}) arise in the context of many physical models, such as the quantum transport and the non-relativistic Newtonian gravity, see \cite{BokaLopezSoler2003,BahGroDonBas2014,Penrose1996} and the references therein.
In particular, if $\lambda=-1$ and $g(x,u)\equiv 0$, Problem (\ref{e:100}) turns out to be a Choquard equation or a nonlinear Schr\"odinger-Newton equation. And if $F(u)=u^2$, Problem (\ref{e:100}) becomes into a Schr\"odinger-Poisson equation. There has been a large number of recent mathematical research on these problems.

For Choquard equation, in the research of Moroz and van Schaftingen \cite{MorozVanSch2013,MorozVanSch2015}, they proved the existence of the ground state solution and verified the related properties. Like the scalar field equation $-\Delta u+u=u^p$, the uniqueness of the positive solution was an open problem for a long time until Ma and Zhao \cite{MaZhao2010} gave a partial affirmative answer. We refer \cite{MorozVanS2017} for a survey. There are also references that focus on more specific solutions. For instance, the existence of multi-bump solutions to a non-autonomous Choquard equation is proved by Wei and Winter \cite{WeiWinter2009}. For the autonomous Choquard equation, \cite{XiaWang2019,WangXia2020} shows the existence of saddle solutions.
Recently, Liu, Ma and Xia \cite{LiuMaXia2021} studied the singularly perturbed Choquard equation and obtained multiple bounded state solutions.

For Schr\"odinger-Poisson equation, Benci and Fortunato \cite{BenciFort1998} studied the corresponding eigenvalue problem. Their attention is focused on Problem (\ref{e:100}) with $g(x,u)=\lambda u$. For problem with nonlinear $g(x,u)$, \cite{Vaira2011} and \cite{IanniRuiz2012} obtained the ground state solution. Schr\"odinger-Poisson equations with critical nonlinear term are studied in \cite{AlvesCassaniTarsiYang2016} for planar case and in \cite{CeramiMolle2019} for the three-dimensional case. There are also numerous references on more specific cases. For example, in \cite{HebeyWei2013}, Hebey and Wei considered the analogue to Schr\"odinger-Poisson equation defined on sphere.

In this paper, we analyze Problem (\ref{e:100}) with $\omega=0$, $F(u)=u$ and $g(x,u)=|u|^{2^* -2}u$, i.e. Problem (\ref{e:000}). The solvability of elliptic problem involving critical exponent was open question for a long time, until the publication of the ground-breaking work by Br\'ezis and Nirenberg \cite{BrezisNirenberg1983}. Subsequent works, such as \cite{CeramiFortStruwe1984,Li1993}, continued to study the multiplicity result for critical elliptic problems and the non-autonomous variants of them. Among all of the problems, we emphasize that the two-dimensional critical problem is with exponentially growing nonlinear terms, cf. \cite{HanLi2010,LiShafrir1994}.

In this paper, we study the critical problem involving potential terms.
The potential term $I_\alpha*_\Omega u$ brings additional difficulties in regularity and the nonlocal eigenvalue problem. To address this problem, we introduce regularity estimates for potential terms. As to the nolocal eigenvalue problem involving term $I_\alpha*_\Omega u$, we apply a minimax approach \cite{RabinowitzBook1986,PereraAgOReBook2010,MolicaRadulaescuSeradeiBook2016} instead of the classical linear functional analytic method \cite[Chapter 6]{EvansBook2010}.

From the point of view of nonlinear analysis, the solution to Problem (\ref{e:000}) is a critical point of the energy functional
\begin{align}
E_{\alpha,\lambda}(u)=\frac{1}{2}\int_\Omega|\nabla u|^2dx-\frac{\lambda}{2}\mathcal{D}_\alpha(u)-\frac{1}{2^*}\int_\Omega |u|^{2^*}dx\nonumber
\end{align}
with
\begin{align}
\mathcal{D}_\alpha(u)=C_{N,\alpha}\int_\Omega dy\int_\Omega \frac{u(x)u(y)}{|x-y|^{N-\alpha}}dy.\nonumber
\end{align}
To deal with $\mathcal{D}_\alpha(u)$, Hardy-Littlewood-Sobolev inequality is an important tool. To be precise, for any $p,q\in(1,\infty)$ and any $\lambda\in(0,n)$, there exists a constant $C_{p,q,\lambda,N}>0$ such that for any $f\in L^p(\mathbb{R}^N)$ and any $g\in L^q(\mathbb{R}^N)$ it holds that
\begin{align}\label{ineq:HLS}
\int_{\mathbb{R}^N}dy \int_{\mathbb{R}^ N}\frac{f(x)g(x)}{|x-y|^\lambda}dx\leq C_{p,q,\lambda,N}|f|_p|g|_q.
\end{align}
This inequality can be traced back to the work by Hardy and Littlewood \cite{HardyLittlewood1928} of 1928. In \cite{Lieb1983}, Lieb studied its sharp constant. We refer \cite{ChenJinLiLim2005,ChenLiOu2005,ChenLi2008,LeiLiMa2012,BellazziniFrankVisciglia2014,DouZhu20151,DouZhu20152} and the references therein for more detailed researches on Inequality (\ref{ineq:HLS}) and its variants. In this paper, we will frequently use Inequality with $p=q=\frac{2N}{N+\alpha}$ and $\lambda=N-\alpha$. The best constant in this case is denoted by $C_{HLS}$.

Now we consider term $\mathcal{D}_\alpha(u)$.
First, notice that $\mathcal{D}_\alpha(u)$ is well-defined on $H_0^1(\Omega)$ since
\begin{align}
\mathcal{D}_\alpha(u) &= C_{N,\alpha}\int_\Omega dy\int_\Omega\frac{u(x)u(y)}{|x-y|^{N-\alpha}}dx=C_{N,\alpha}\int_{\mathbb{R}^N} dy\int_{\mathbb{R}^N}\frac{\widetilde{u}(x)\widetilde{u}(y)}{|x-y|^{N-\alpha}}dx\nonumber\\
&\leq C|\widetilde{u}|_{\frac{2N}{N+\alpha}}^2\leq C|\nabla u|_2^2.\nonumber
\end{align}
Here, $\widetilde{u}$ is the zero extension of $u$, i.e.
\begin{equation}
\widetilde{u}(x)=\left\{
\begin{aligned}
u(x) & \qquad & x\in\Omega, \nonumber\\
0 & \qquad & x\notin\Omega. \nonumber
\end{aligned}
\right.
\end{equation}
By a similar approach, we can prove that $\mathcal{D}_\alpha(u)$ is a $C^2$ functional.
\subsection{Main results}
The first part of our result concerns the non-existence of solutions to Problem (\ref{e:000}).
\begin{theorem}\label{t:nonexistence}
Suppose that $\Omega$ is star shaped and that $\lambda\leq0$, then Problem (\ref{e:000}) admits no non-trivial solutions.
\end{theorem}
Under certain assumptions, there is also a nonexistence result for Problem (\ref{e:000}) for some positive $\lambda>0$.
\begin{theorem}\label{t:nonexistence2}
Suppose that $N-\alpha-4<0$ and that there exists a constant $c_0>0$ such that $x\cdot n_x\geq c_0$ for any $x\in\partial\Omega$, then there exists a constant $\Lambda_0=\Lambda_0(\alpha,\Omega)>0$ such that for any $\lambda\leq\Lambda_0$, Problem (\ref{e:000}) admits no positive solution.
\end{theorem}

On the other hand, we establish the existence result.
\begin{theorem}\label{t:existence}
Problem (\ref{e:000}) admits a solution if one of the following cases occurs.
\begin{itemize}
  \item [$(1).$] $N-\alpha-4\geq0$ and $\lambda>0$;
  \item [$(2).$] $N-\alpha-4<0$ and $\lambda>0$ large.
\end{itemize}
If, additionally, we assume that $0<\lambda<\lambda_1$ with
\begin{align}
\lambda_1=\inf_{u\neq0}\frac{\int_\Omega|\nabla u|^2 dx}{\mathcal{D}_\alpha(u)},\nonumber
\end{align}
Problem (\ref{e:000}) admits a positive solution.
\end{theorem}

To discuss the multiplicity results, we need the following nonlocal eigenvalue problem.
\begin{equation}\label{e:eigen}
    \left\{
   \begin{array}{lr}
     -{\Delta}u=\lambda I_\alpha*_\Omega u\mbox{ in }\Omega ,\\
     \lambda\in\mathbb{R},\,u\in H_0^1(\Omega).
   \end{array}
   \right.
\end{equation}
For Problem (\ref{e:eigen}), the following theorem holds.
\begin{theorem}\label{t:eigen}
There exists a sequence $\{(\lambda_k,\phi_k)\}_k\subset \mathbb{R}\times H_0^1(\Omega)$ such that
\begin{itemize}
  \item [$(1).$] $(\lambda_k,\phi_k)$ solves Problem (\ref{e:eigen});
  \item [$(2).$] $\phi_k\in C^{2,\theta}(\Omega)\cap C_0(\overline{\Omega})$ for some $\theta\in(0,1)$ and $\{\phi_k\}_k$ forms an normalized orthogonal base of $H_0^1(\Omega)$;
  \item [$(3).$] $0<\lambda_1<\lambda_2\leq\lambda_3\leq\cdots\to\infty$ and $\phi_1>0$.
\end{itemize}
\end{theorem}

Based on Theorem \ref{t:eigen}, we get the following multiplicity results.
\begin{theorem}\label{t:multiplicity}
Under the assumption of Theorem \ref{t:existence}, denote
\begin{itemize}
  \item $\upsilon=S\cdot|\Omega|_N^{-\frac{2+\alpha}{N}}\cdot C_{HLS}^{-1}C_{N,\alpha}^{-1}$;
  \item $m=\#\{j\in\mathbb{Z}_+|\lambda<\lambda_j<\lambda+\upsilon\}$.
\end{itemize}
Then, Problem (\ref{e:000}) admits at least $m$ pairs of solutions. Here,
$S=\inf_{u\neq0}\frac{\int_{\mathbb{R}^N}|\nabla u|^2 dx}{(\int_{\mathbb{R}^N}|u|^{2^*}dx)^\frac{2}{2^*}}$
and $\lambda_k$ are the eigenvalues in Theorem \ref{e:eigen}.
\end{theorem}

\subsection{Organization of this paper}
In Section \ref{Sec:Eigenproblem}, we prove Theorem \ref{t:eigen}. Theorems \ref{t:nonexistence} and Theorem \ref{t:nonexistence2} are proved in Section \ref{Sec:Nonexistence}. To this end, we first show the regularity of the weak solutions to Problem (\ref{e:000}) in Subsection \ref{Subsec:regular} and in Subsection \ref{Subsec:Pohozaev} we prove the corresponding Pohozaev identity. In Section \ref{Sec:Existence}, Theorem \ref{t:existence} is proved by testing the mountain pass level and verifying a local compactness result and a maximum principle for a nonlocal operator. In order to prove \ref{t:multiplicity}, we need to estimate the genus of certain energy level sets. This is shown in Subsection \ref{Subsec:multiple}.

\vs

\noindent{\bf Notations.}
\begin{itemize}
  \item Throughout this paper, $C$ denotes a generic positive constant that may vary from line to line;
  \item $o_n(1)$ denotes a generic infinitesimal as $n\to\infty$;
  \item $dS_x$ denotes the area element of $\partial\Omega$ and $d\vec{S}_x=n_x dS_x$ where $n_x$ is the outer normal vector of $\partial\Omega$ at $x$;
  \item Without causing confusion, the norm of the spaces $L^p(\Omega)$ and of the space $H_0^1(\Omega)$ is denoted by $|\cdot|_p$ and $\|\cdot\|$, respectively;
  \item The constant $C_{N,\alpha}=\frac{\Gamma(\frac{N-\alpha}{2})}{\Gamma(\frac{\alpha}{2})\pi^\frac{N}{2}2^\alpha}$ and $C_{HLS}$ denotes the best constant of Inequality (\ref{ineq:HLS}) with $p=q=\frac{2N}{N+\alpha}$ and $\lambda=N-\alpha$;
  \item $S=\inf\{\int_{\mathbb{R}^N}|\nabla u|^2dx|\int_{\mathbb{R}^N}|u|^{2^*}dx=1\}>0$.
\end{itemize}

\section{Proof of Theorem \ref{t:eigen}: The eigenvalue problems}\label{Sec:Eigenproblem}
In this part, we prove Theorem \ref{t:eigen} in steps. Different from the classical approach \cite[Chapter 6]{EvansBook2010}, we apply a minimax approach which can also be applied to nonlinear eigenvalue problems and other nonlocal eigenvalue problems, cf. \cite{RabinowitzBook1986,PereraAgOReBook2010,MolicaRadulaescuSeradeiBook2016}.

\subsection{The existence of eigenfunctions}
\begin{lemma}
The set $\mathcal{M}=\{u\in H_0^1(\Omega)|\mathcal{D}_\alpha(u)=1\}$ is a $C^1$ manifold in $H_0^1(\Omega)$ that is homomorphic to $S$. Here $S=\{u\in H_0^1(\Omega)|\int_\Omega|\nabla u|^2 dx=1\}$.
\end{lemma}
\noindent{\bf Proof.}
To verify this, it is sufficient to notice that
\begin{align}
\mathcal{D}_\alpha(u)=C_{N,\alpha} \int_\Omega dy\int_\Omega\frac{u(x)u(y)dx}{|x-y|^{N-\alpha}}= C_{N,\alpha} \int_{\mathbb{R}^N} dy\int_{\mathbb{R}^N}\frac{\widetilde{u}(x)\widetilde{u}(y)dx}{|x-y|^{N-\alpha}}.\nonumber
\end{align}
Here, $\widetilde{u}$ is the zero-extension of $u$ outside of $\Omega$. Notice that
\begin{align}
|\mathcal{D}_\alpha(u)|\leq C|\widetilde{u}|^2_\frac{2N}{N+\alpha}\leq C|\nabla u|_2^2.\nonumber
\end{align}
This is due to Hardy-Littlewood-Sobolev inequality and H\"older's inequality. This implies that $\mathcal{D}_\alpha (u)<\infty$ for any $H_0^1(\Omega)$. On the other hand,
by Parseval's identity,
\begin{align}
\mathcal{D}_\alpha(u)=C_{N,\alpha}\int_{\mathbb{R}^N}|I_\frac{\alpha}{2}*\widetilde{u}|^2 dx.\nonumber
\end{align}
Therefore, $D_\alpha(u)=0$ if and only if $\widetilde{u}=0$, i.e. $u=0$ in $\Omega$.
\begin{flushright}
$\Box$
\end{flushright}

Another necessary result is that
\begin{lemma}
$\int_\Omega|\nabla u|^2dx$ satisfies $(PS)$ condition on $\mathcal{M}$. To be precise, if the sequence $\{u_n\}_n\subset H_0^1(\Omega)$ and $\{\lambda_n\}_n\subset\mathbb{R}$ and a number $\lambda_*\in\mathbb{R}$ satisfy
\begin{itemize}
  \item [$(a).$] $\int_\Omega|\nabla u_n|^2dx\to\lambda_*$;
  \item [$(b).$] $u_n-\lambda_n(-\Delta)_\Omega^{-1}I_\alpha*_\Omega u_n=o_n(1)$,
\end{itemize}
then $\lambda_n\to\lambda_*$ and there exists $u_\infty\in \mathcal{M}$ such that $u_n\to u_\infty$ in $H_0^1(\Omega)$.
\end{lemma}
\noindent{\bf Proof.}
Notice that $\{u_n\}_n$ is bounded in $H_0^1(\Omega)$, then define $u_\infty$ as the weak limit. By a direct computation,
\begin{align}\label{convergence1}
\mathcal{D}_\alpha(u_n -u_\infty)\leq C|u_n-u_\infty|_\frac{2N}{N+\alpha}\to 0
\end{align}
due to the Sobolev embedding and
\begin{align}
\mathcal{D}_\alpha(u_n -u_\infty)=\mathcal{D}_\alpha(u_n)+\mathcal{D}_\alpha(u_\infty)-2C_{N,\alpha}\int_\Omega dy \int_\Omega\frac{u_n(x)u_\infty(y)}{|x-y|^{N-\alpha}}dx.\nonumber
\end{align}
Since
\begin{itemize}
  \item $u_n\to u_\infty$ in $L^2(\Omega)$;
  \item $I_\alpha*_\Omega u_\infty=I_\alpha*_{\mathbb{R}^N}\widetilde{u}_\infty|_\Omega\in L^2(\Omega)$. Here, $\widetilde{u}_\infty$ is the zero extension of $u_\infty$,
\end{itemize}
we get
\begin{align}\label{convergence2}
\mathcal{D}_\alpha(u_n -u_\infty)=\mathcal{D}_\alpha(u_n)-\mathcal{D}_\alpha(u_\infty).
\end{align}
Moreover,
\begin{align}\label{convergence3}
\mathcal{D}_\alpha(u_n -u_\infty)\geq0.
\end{align}
Then, (\ref{convergence1}), (\ref{convergence2}) and (\ref{convergence3}) imply that $\mathcal{D}_\alpha(u_\infty)=1$.
Multiplying $(b)$ by $u_n$, we get
\begin{align}
\int_\Omega|\nabla u_n|^2dx=\lambda_n\mathcal{D}_\alpha(u_n)+o_n(1).\nonumber
\end{align}
Combining $(a)$, we get $\lambda_n\to\lambda_*$. Multiplying $(b)$ by $u_\infty$,
\begin{align}
\int_\Omega\nabla u_n\cdot \nabla u_\infty dx-\lambda_n C_{N,\alpha}\int_\Omega dy\int_\Omega\frac{u_n(x)u_\infty( y)}{|x-y|^{N-\alpha}}dx+o_n(1)=0.\nonumber
\end{align}
Then,
\begin{align}
\|u_n-u_\infty\|^2=\int_\Omega|\nabla u_n|^2dx-\int_\Omega|\nabla u_\infty|^2dx=\lambda_n -\lambda_*+o_n(1)\to0.\nonumber
\end{align}
This completes the proof.

\begin{flushright}
$\Box$
\end{flushright}

Consider the following minimax problem
\begin{align}\label{minimax:eigen}
\lambda_k=\inf_{A\subset\mathcal{M},\gamma(A)\geq k}\sup_{u\in A}\frac{\int|\nabla u|^2dx}{\mathcal{D}_\alpha(u)}.
\end{align}
Here, $\gamma$ is the $\mathbb{Z}_2$-genus generated by the anti-podal action $u\mapsto-u$. By a standard method (see, for instance, \cite{AmbrosettiRabinowitz1973} and \cite{RabinowitzBook1986}), we obtain a sequence $\{(\lambda_k,\phi_k)\}_k\subset\mathbb{R}\times\mathcal{M}$ solving Problem (\ref{e:eigen}) and $\lambda_k\to+\infty$.

Especially,
\begin{align}
\lambda_1=\inf_{u\in\mathcal{M}}\frac{\int|\nabla u|^2dx}{\mathcal{D}_\alpha(u)}.\nonumber
\end{align}
By Hardy-Littlewood-Sobolev inequality,
\begin{align}
\mathcal{D}_\alpha(u)&=C_{N,\alpha}\int_\Omega dy\int_\Omega\frac{u(x)u(y)}{|x-y|^{N-\alpha}}dx =C_{N,\alpha}\int_{\mathbb{R}^N} dy\int_{\mathbb{R}^N} \frac{\widetilde{u}(x)\widetilde{u}(y)}{|x-y|^{N-\alpha}}dx\nonumber\\
&\leq C|\widetilde{u}|^2_{L^\frac{2N}{N+\alpha}(\mathbb{R}^N)}\leq C|\nabla u|_2^2.\nonumber
\end{align}
Here, $\widetilde{u}$ is the zero-extension of $u$. \cite[Theorem 8.10]{RabinowitzBook1986} and \cite{AmbrosettiRabinowitz1973} imply the following result.
\begin{proposition}\label{prop:eigenfunctions}
There exists a sequence of solutions $\{(\lambda_k,\phi_k)\}_k$ to Problem (\ref{e:eigen}) such that  $0<\lambda_1\leq\lambda_2\leq\lambda_3\leq\cdots$.
\end{proposition}

\subsection{The regularity of eigenfunctions}
In this part, we prove a regularity result for the solutions to Problem (\ref{e:eigen}). We start with a function $u\in H_0^1(\Omega)$ such that
\begin{align}\label{equ:eigen}
-\Delta u=\lambda C_{N,\alpha}\int_\Omega\frac{u(y)}{|x-y|^{N-\alpha}}dy
\end{align}
for some $\lambda>0$. Recall \cite[Lemma 7.12]{GilbargTrudinger2011}.
\begin{lemma}\label{l:GT}
The operator $I_\alpha*_\Omega \cdot:L^p(\Omega)\to L^q(\Omega)$ continuously for any $q\in[1,\infty]$ and $p$ satisfying $0\leq\frac{1}{p}-\frac{1}{q}<\frac{\alpha}{N}$. To be precise, there exists a constant $C>0$ such that $|I_\alpha*_\Omega u|_q\leq C|u|_p$ for any $u\in L^p(\Omega)$.
\end{lemma}
In follows, it is sufficient to use the case $p=q$.
\begin{lemma}\label{l:RegularLp}
Let $(\lambda,u)$ solving Problem (\ref{equ:eigen}).
For any $r\geq1$, we get $u\in W^{2,r}(\Omega)$.
\end{lemma}
\noindent{\bf Proof.}
To begin with, let us recall that $u\in H_0^1(\Omega)$. Then, $u\in L^\frac{2N}{N-2}(\Omega)$ and Sobolev inequality and Lemma \ref{l:GT} give

\vs

\begin{itemize}
  \item [$\mbox{Step 1}.$] $u\in W^{2,\frac{2N}{N-2}}(\Omega)$, then $u\in L^\frac{2N}{N-6}(\Omega)$ and $I_\alpha*_\Omega u\in L^\frac{2N}{N-6}(\Omega)$;
  \item [$\mbox{Step 2}.$] $u\in W^{2,\frac{2N}{N-6}}(\Omega)$, then $u\in L^\frac{2N}{N-10}$ and hence $I_\alpha*_\Omega u\in L^\frac{2N}{N-10}(\Omega)$;
  \item [$ $] ......
  \item [$\mbox{Step k}.$] $u\in W^{2,\frac{2N}{N-4k+2}}(\Omega)$, then $u\in L^\frac{2N}{N-4k-2}(\Omega)$ and hence $I_\alpha*_\Omega u\in L^\frac{2N}{N-4k-2}(\Omega)$;
  \item [$ $] ......
\end{itemize}

\vs

Let $k_m:=\max\{k\in \mathbb{N}|N-4k>2\}$, then

\vs

\begin{itemize}
  \item [$\mbox{Case 1}.$] $2=N-4(k_m+1)$, then for any $r\geq 1$, $u\in L^r(\Omega)$;
  \item [$\mbox{Case 2}.$] $2>N-4(k_m+1)$, then $u\in L^\infty(\Omega)$ and for any $r\geq 1$, $u\in L^r(\Omega)$.
\end{itemize}

\vs

This implies that $u\in W^{2,r}(\Omega)$ for any $r\geq1$.

\begin{flushright}
$\Box$
\end{flushright}

To proceed, we need to introduce the $M^p(\Omega)$ space and the associated potential estimate.
\begin{definition}
$M^p(\Omega)=\{u\in L^1_{loc}(\Omega)|\exists K>0\mbox{ s.t. }\forall x\in\Omega\mbox{ and }R>0,\, \int_{\Omega\cap B_R(x)}|u|dx\leq K R^{N(1-\frac{1}{p})}\}$. The norm of $M^p(\Omega)$ is defined as
\begin{align}
\|u\|_{M^p(\Omega)}=\inf\Big\{K>0\Big|\forall x\in\Omega\mbox{ and }R>0,\, \int_{\Omega\cap B_R(x)}|u|dx\leq K R^{N(1-\frac{1}{p})}\Big\}.\nonumber
\end{align}
\end{definition}
\begin{remark}
It holds that $L^1(\Omega)=M^1(\Omega)$, $L^p(\Omega)\subset M^p(\Omega)$ and $L^\infty(\Omega)=M^\infty(\Omega)$. We refer \cite[Chapter 7.9]{GilbargTrudinger2011}.
\end{remark}
The following result is \cite[Lemma 7.18]{GilbargTrudinger2011}.
\begin{lemma}\label{l:potential2}
If $\frac{N}{\alpha}<p$, there exists a constant $C>0$, $|I_\alpha*_\Omega u|_\infty\leq C\|u\|_{M^p(\Omega)}$ for any $u\in M^p(\Omega)$.
\end{lemma}

\begin{lemma}\label{l:RegularGradient}
Let $(\lambda,u)$ solving Problem (\ref{equ:eigen}).
Then, we get $\nabla I_\alpha*_\Omega u\in L^\infty(\Omega)$.
\end{lemma}
\noindent{\bf Proof.}
Using Lemma \ref{l:RegularLp}, we get a $\theta\in(0,1)$ such that $u\in C_0^{1,\theta}(\Omega)$. Then, $|\nabla u|\in C^{0,\theta}(\Omega)\subset L^\infty(\Omega)$. On the other hand, for any $x\in\Omega$,
\begin{align}
\nabla_x\int_\Omega\frac{u(y)}{|x-y|^{N-\alpha}}dy &= \nabla_x\int_{\mathbb{R}^N} \frac{\widetilde{u}(y)}{|x-y|^{N-\alpha}}dy= \int_{\mathbb{R}^N} \frac{\nabla\widetilde{u}(y)}{|x-y|^{N-\alpha}}dy\nonumber\\
&=\int_{\Omega} \frac{\nabla u(y)}{|x-y|^{N-\alpha}}dy.\nonumber
\end{align}
Here, $\widetilde{u}$ is the zero extension of $u$.
Then, $|\nabla_x\int_\Omega\frac{u(y)}{|x-y|^{N-\alpha}}dy|\leq \int_{\Omega} \frac{|\nabla u(y)|}{|x-y|^{N-\alpha}}$ for any $x\in\Omega$. Then, Lemma \ref{l:potential2} gives that $\nabla I_\alpha*_\Omega u\in L^\infty(\Omega)$.

\begin{flushright}
$\Box$
\end{flushright}

\begin{proposition}\label{prop:eigenregular}
For any solution $(\lambda,u)$ to Problem (\ref{equ:eigen}), $u\in C^{2,\theta}(\Omega)\cap C_0(\overline{\Omega})$ for some $\theta\in(0,1)$.
\end{proposition}
\noindent{\bf Proof.}
Lemma \ref{l:RegularGradient} implies that $I_\alpha*_\Omega u\in C^{0,\theta}(\Omega)$ for some $\theta\in(0,1)$. The lemma follows the classical H\"older estimate \cite[Chapter 6]{GilbargTrudinger2011}.

\begin{flushright}
$\Box$
\end{flushright}

\subsection{Orthogonality and completeness}
The following result is evident.
\begin{claim}
For the sequence we obtain in Proposition \ref{prop:eigenfunctions}, it holds that
\begin{itemize}
  \item [$(1).$] For any $k$ and any $\mu\in\mathbb{R}$, $(\lambda_k,\mu \phi_k)$ solves Problem (\ref{e:eigen});
  \item [$(2).$] If there exists $k_0,p_0$ such that $\lambda_{k_0+1}=\cdots=\lambda_{k_0+p_0}$, then for any $\mu_1,\cdots,\mu_{p_0}$, $(\lambda_{k_0+1},\sum_{i=1}^{p_0}\mu_i\phi_i)$ solves Problem (\ref{e:eigen}).
\end{itemize}
\end{claim}

Furthermore, we get
\begin{lemma}\label{l:orthogonal}
If $i,j=1,2,\cdots$ such that $\lambda_i\neq\lambda_j$, then $\int_\Omega\nabla \phi_i\cdot \nabla \phi_j dx=0$.
\end{lemma}
\noindent{\bf Proof.}
Notice that
\begin{align}\label{equ:i}
-{\Delta}\phi_i=\lambda_i I_\alpha*_\Omega \phi_i
\end{align}
and
\begin{align}\label{equ:j}
-\Delta \phi_j=\lambda_j I_\alpha*_\Omega \phi_j.
\end{align}
Multiplying (\ref{equ:i}) by $\frac{1}{\lambda_i}\phi_j$, (\ref{equ:j}) by $\frac{1}{\lambda_j}\phi_i$, integrating over $\Omega$ and taking the difference, we get
\begin{align}
\int_\Omega\nabla \phi_i\cdot \nabla \phi_j dx=0.\nonumber
\end{align}
\begin{flushright}
$\Box$
\end{flushright}

\begin{remark}
Moreover, under the assumptions of Lemma \ref{l:orthogonal}, we also get $\int_\Omega dy\int_\Omega\frac{\phi_i(x)\phi_j(y)}{|x-y|^{N-\alpha}}dx=0$.
\end{remark}

Therefore, without loss of generality, we can assume that $\{\phi_k\}_k$ are pairwise orthogonal. By the $(PS)$ condition, it is also clear that
\begin{claim}
For any $k=1,2,\cdots$, $\sup\{l>0|\lambda_{k+1}=\cdots=\lambda_{k+l}\}<\infty$.
\end{claim}

Now we verify the completeness of $\{\phi_k\}_k$.
\begin{proposition}\label{prop:complete}
It holds that $H_0^1(\Omega)=\overline{\mbox{span}\{\phi_1,\cdots,\phi_k,\cdots\}}$.
\end{proposition}
\noindent{\bf Proof.}
Otherwise, if there exists a function $\phi_*\in H_0^1(\Omega)$ such that $\phi_*\perp\overline{\mbox{span}\{\phi_1,\cdots,\phi_k,\cdots\}}$.
Denote $H_0:=\overline{\mbox{span}\{\phi_1,\cdots,\phi_k,\cdots\}}$ and $H_0^\perp:=(\overline{\mbox{span}\{\phi_1,\cdots,\phi_k,\cdots\}})^\perp$. Then $\phi_*\in H_0^\perp$.
Now we claim that
\begin{claim}
$\mathcal{N}=\mathcal{M}\cap H_0^\perp$ is a natural constrain for $\int|\nabla u|^2 dx$ on $\mathcal{M}$.
\end{claim}

If $\nabla_{\mathcal{N}}\int|\nabla u|^2dx=0$ for $u\in\mathcal{N}$, we want to prove $\nabla_{\mathcal{M}}\int|\nabla u|^2dx=0$. To do this, let us select a $\phi_j\in H_0$. Then,
\begin{align}
0=\int_\Omega\nabla u\cdot\nabla \phi_j dx=\int_\Omega u(-\Delta \phi_j) dx=\lambda_j C_{N,\alpha}\int_\Omega dy\int_\Omega\frac{u(x)\phi_j(y)}{|x-y|^{N-\alpha}}dx,\nonumber
\end{align}
since $u\perp\phi_j$ in $H_0^1(\Omega)$ and $\phi_j$ is the $j$-th eigenfunction. In this computation, we also proved that $v\in T_u\mathcal{M}$. The same computation implies that $\nabla_\mathcal{M}\int|\nabla u|^2dx=0$ and the claim is proved.

Then, we obtain an eigenvalue and an eigenfunction $(\lambda_0,u_0)$ with $0<\lambda_0<\infty$. This contradicts with the order of the eigenvalues $\lambda_1\leq\lambda_2\leq\cdots$. Therefore, $H_0^\perp=0$, i.e. $H_0^1(\Omega)=\overline{\mbox{span}\{\phi_1,\cdots,\phi_k,\cdots\}}$.

\begin{flushright}
$\Box$
\end{flushright}

\subsection{Principal eigenfunction}
In this part, we analyze the principal eigenfunction and the principle eigenvalue.
Recall that
\begin{align}\label{eign:principal}
\lambda_1=\inf_{u\in\mathcal{M}}\int_\Omega|\nabla u|^2dx
\end{align}
and $\phi_1$ attends the above minimum.

\begin{proposition}\label{prop:phi1>0}
It holds that $\phi_1>0$.
\end{proposition}
\noindent{\bf Proof.}
Otherwise, we assume that $\phi_1$ changes its sign. Let $\phi_{1,+}=\max\{\phi_1,0\}$ and $\phi_{1,-}=\min\{\phi_1,0\}$.
Then,
\begin{align}
\int_\Omega|\nabla \phi_1|^2dx&=\lambda_1 C_{N,\alpha}\int_\Omega dy\int_\Omega\frac{\phi_1(x) \phi_1(y)}{|x-y|^{N-\alpha}}dx\nonumber\\
&=\lambda_1 C_{N,\alpha}\int_\Omega dy\int_\Omega\frac{\phi_{1,+}(x) \phi_{1,+}(y)}{|x-y|^{N-\alpha}}dx
+\lambda_1 C_{N,\alpha}\int_\Omega dy\int_\Omega\frac{\phi_{1,-}(x) \phi_{1,-}(y)}{|x-y|^{N-\alpha}}dx\nonumber\\
&\quad +2\lambda_1 C_{N,\alpha}\int_\Omega dy\int_\Omega\frac{\phi_{1,+}(x) \phi_{1,-}(y)}{|x-y|^{N-\alpha}}dx.\nonumber
\end{align}
Here,
\begin{align}
\lambda_1 C_{N,\alpha}\int_\Omega dy\int_\Omega\frac{\phi_{1,+}(x) \phi_{1,-}(y)}{|x-y|^{N-\alpha}}dx\leq0.\nonumber
\end{align}
If $\lambda_1 C_{N,\alpha}\int_\Omega dy\int_\Omega\frac{\phi_{1,+}(x) \phi_{1,-}(y)}{|x-y|^{N-\alpha}}dx<0$, then
\begin{align}
\int_\Omega|\nabla \phi_1|^2dx<\lambda_1 C_{N,\alpha}\int_\Omega dy\int_\Omega\frac{\phi_{1,+}(x) \phi_{1,-}(y)}{|x-y|^{N-\alpha}}dx
+\lambda_1 C_{N,\alpha}\int_\Omega dy\int_\Omega\frac{\phi_{1,-}(x) \phi_{1,-}(y)}{|x-y|^{N-\alpha}}dx.\nonumber
\end{align}
Then, either
\begin{align}
\int_\Omega|\nabla \phi_{1,+}|^2dx<\lambda_1 C_{N,\alpha}\int_\Omega dy\int_\Omega\frac{\phi_{1,+}(x) \phi_{1,+}(y)}{|x-y|^{N-\alpha}}dx\nonumber
\end{align}
or
\begin{align}
\int_\Omega|\nabla \phi_{1,-}|^2dx<\lambda_1 C_{N,\alpha}\int_\Omega dy\int_\Omega\frac{\phi_{1,-}(x) \phi_{1,-}(y)}{|x-y|^{N-\alpha}}dx.\nonumber
\end{align}
This contradict with (\ref{eign:principal}). Then,
\begin{align}
\lambda_1 C_{N,\alpha}\int_\Omega dy\int_\Omega\frac{\phi_{1,+}(x) \phi_{1,-}(y)}{|x-y|^{N-\alpha}}dx=0.\nonumber
\end{align}
Hence,
\begin{align}
\int_\Omega|\nabla \phi_1|^2dx=\lambda_1 C_{N,\alpha}\int_\Omega dy\int_\Omega\frac{\phi_{1,+}(x) \phi_{1,-}(y)}{|x-y|^{N-\alpha}}dx
+\lambda_1 C_{N,\alpha}\int_\Omega dy\int_\Omega\frac{\phi_{1,-}(x) \phi_{1,-}(y)}{|x-y|^{N-\alpha}}dx.\nonumber
\end{align}
Then, $\phi_{1,+}$ and $\phi_{1,-}$ also minimize (\ref{eign:principal}) and
\begin{align}
-\Delta \phi_{1,\pm}=\lambda_1 I_\alpha*_\Omega\phi_{1,\pm}\mbox{ in }\Omega.\nonumber
\end{align}
This contradicts with the strong maximum principle of $-\Delta$, cf. \cite[Chapter 6.4]{EvansBook2010}.

\begin{flushright}
$\Box$
\end{flushright}

\subsection{Proof of Theorem \ref{t:eigen}}
The proof of Theorem \ref{t:eigen} follows the above considerations.

\vs

\noindent{\bf Proof of Theorem \ref{t:eigen}.}
Assertion $(1)$ of Theorem \ref{e:eigen} is proved by Proposition \ref{prop:eigenfunctions}, i.e. there exists a sequence of pairs $\{(\lambda_k,\phi_k)\}_k$ whose elements solve Problem (\ref{e:eigen}). By Proposition \ref{prop:eigenregular}, for any $k=1,2,\cdots$, $\phi_k\in C^{2,\theta}(\Omega)\cap C_0(\overline{\Omega})$ for some $\theta\in(0,1)$. For the set $\{\phi_k\}_k$, the orthogonality and the completeness are guaranteed by Lemma \ref{l:orthogonal} and Proposition \ref{prop:complete}. This proves Assertion $(2)$ of Theorem \ref{e:eigen}.

As for Assertion $(3)$, using Proposition \ref{prop:phi1>0}, the principal eigenvalue $\lambda_1$ is simple. The proof is complete.

\begin{flushright}
$\Box$
\end{flushright}

\section{Proof of Theorems \ref{t:nonexistence} and \ref{t:nonexistence2}: Non-existence results}\label{Sec:Nonexistence}
\subsection{The regularity of weak solutions}\label{Subsec:regular}
The proofs of Theorems \ref{t:nonexistence} and \ref{t:nonexistence2} are based on \cite{BrezisNirenberg1983}. We start with $u\in H^1_0(\Omega)$ solving Problem (\ref{e:000}), i.e.,
\begin{equation}
    \left\{
   \begin{array}{lr}
     -{\Delta}u=\lambda I_\alpha*_\Omega u+|u|^{2^*-2}u\mbox{ in }\Omega ,\nonumber\\
     u\in H_0^1(\Omega).\nonumber
   \end{array}
   \right.
\end{equation}
We want to prove that $u\in C^{1,\theta}(\Omega)$ for some $\theta\in(0,1)$. Recall the following Hardy-Littlewood-Sobolev inequality:
\begin{align}
\Big|\int_{\mathbb{R}^N}dy\int_{\mathbb{R}^N}\frac{f(x)g(y)}{|x-y|^{N-\alpha}}dx\Big|\leq C|f|_p|g|_t\nonumber
\end{align}
with $\frac{1}{p}+\frac{1}{t}=1+\frac{\alpha}{N}$. First, we have that
\begin{lemma}\label{l:nonlinearregular1}
Let $u$ be a solution to Problem (\ref{e:000}) in $H_0^1(\Omega)$.
For any $r\geq 1$, $u\in L^r(\Omega)$.
\end{lemma}
\noindent{\bf Proof.}
Multiplying Problem (\ref{e:000}) by $|u|^{s-1}u$ and integrating over $\Omega$, we get
\begin{align}
\frac{4s}{(s+1)^2}\int_\Omega|\nabla (u^\frac{s+1}{2})|^2dx=\int_\Omega(-\Delta u)|u|^{s-1}u dx\leq \lambda C\int_{\mathbb{R}^N}dy
\int_{\mathbb{R}^N}\frac{|\widetilde{u}(y)\widetilde{u}(x)^s|}{|x-y|^{N-\alpha}}dx+ |u|_{\frac{N+2}{N-2}+s}^{\frac{N+2}{N-2}+s}.\nonumber
\end{align}
Here, $\widetilde{u}$ is the zero extension of $u$.
Using Sobolev inequality and Hardy-Littlewood-Sobolev inequality,
\begin{align}
|u|_\frac{(s+1)N}{N-2}^s\leq C\Big(|u|^{s+1}_\frac{N(s+1)}{N+\alpha}+ |u|_{\frac{N+2}{N-2}+s}^{\frac{N+2}{N-2}+s}\Big).\nonumber
\end{align}
Noticing $\frac{N(s+1)}{N+\alpha}\leq \frac{N+2}{N-2} +s$, then
\begin{align}\label{ineq:iteration}
|u|_\frac{(s+1)N}{N-2}^s\leq C\Big(|u|^{s+1}_{\frac{N+2}{N-2}+s}+ |u|_{\frac{N+2}{N-2}+s}^{\frac{N+2}{N-2}+s}\Big).
\end{align}
Letting $s_{i+1}=\frac{N}{N-2}\cdot s_i-\frac{2}{N-2}$. We can prove that for any $r\geq1$, $u\in L^r(\Omega)$.
\begin{flushright}
$\Box$
\end{flushright}

\begin{lemma}\label{l:RegularNonlinearC2a}
Let $u$ be a solution to Problem (\ref{e:000}) in $H_0^1(\Omega)$.
There exists a $\theta\in(0,1)$, $u\in C^{2,\theta}(\Omega)$.
\end{lemma}
\noindent{\bf Proof.}
By Lemma \ref{l:nonlinearregular1} and \ref{l:GT}, we get for any $r\geq 1$, $u\in W^{2,r}(\Omega)$. Then, for any $r\geq 1$, $|\nabla u|\in L^r(\Omega)$. Notice that for any $x\in\Omega$,
\begin{align}
\nabla_x\int_\Omega\frac{u(y)}{|x-y|^{N-\alpha}}dy &= \nabla_x\int_{\mathbb{R}^N} \frac{\widetilde{u}(y)}{|x-y|^{N-\alpha}}dy= \int_{\mathbb{R}^N} \frac{\nabla\widetilde{u}(y)}{|x-y|^{N-\alpha}}dy\nonumber\\
&=\int_{\Omega} \frac{\nabla u(y)}{|x-y|^{N-\alpha}}dy.\nonumber
\end{align}
Lemma \ref{l:RegularGradient} implies that $|\nabla_x \int_\Omega\frac{u(y)}{|x-y|^{N-\alpha}}dy|\in L^\infty(\Omega)$. Then, $I_\alpha*_\Omega u\in C^{0,\theta}(\Omega)$ for some $\theta\in(0,1)$. Combing the classical H\"older estimate (cf. \cite[Chapter 6]{GilbargTrudinger2011}), we complete the proof.

\begin{flushright}
$\Box$
\end{flushright}

\subsection{Pohozaev identity}\label{Subsec:Pohozaev}
In this part, we derive a Pohozaev identity for Problem (\ref{e:000}). First, let us consider the potential term $I_\alpha*_\Omega u$.
\begin{lemma}\label{l:PotentialPohozaev}
It holds that
\begin{align}
\int_\Omega(I_\alpha*_\Omega u)(x\cdot\nabla u)dx =-\frac{N+\alpha}{2} C_{N,\alpha} \int_\Omega dx\int_\Omega\frac{u(x)u(y)}{|x-y|^{N-\alpha}}dy.\nonumber
\end{align}
\end{lemma}
\noindent{\bf Proof.}
By a direct computation,
\begin{align}
\mbox{div}_x\Big((I_\alpha*_\Omega u)x\cdot u\Big) =\mbox{div}_x\Big((I_\alpha*_\Omega)x\Big)u+(I_\alpha*_\Omega)x\cdot\nabla u.\nonumber
\end{align}
This gives that
\begin{align}
\int_\Omega(I_\alpha*_\Omega u)(x\cdot\nabla u)dx &=
-\int_\Omega dx\Big[\mbox{div}_x\Big((I_\alpha*_\Omega u)x\Big)u\Big]+\int_\Omega dx\Big[\mbox{div}_x\Big((I_\alpha*_\Omega u)xu(x)\Big)\Big]\nonumber\\
&= -\int_\Omega dx\Big[\mbox{div}_x\Big((I_\alpha*_\Omega u)x\Big)u\Big]\nonumber
\end{align}
since $u\in H_0^1(\Omega)$. Moreover, notice that
\begin{align}
\mbox{div}_x\Big(\frac{x}{|x-y|^{N-\alpha}}\Big)=\frac{N}{|x-y|^{N-\alpha}}-(N-\alpha)\frac{(x-y)\cdot x}{|x-y|^{N-\alpha+2}}.\nonumber
\end{align}
We get
\begin{align}
\int_\Omega(I_\alpha*_\Omega u)(x\cdot\nabla u)dx & = -N C_{N,\alpha}\int_\Omega dx\int_\Omega\frac{u(x)u(y)}{|x-y|^{N-\alpha}}dy\nonumber\\
&\qquad+(N-\alpha)C_{N,\alpha}\int_\Omega dx\int_\Omega\frac{(x-y)\cdot xu(x)u(y)}{|x-y|^{N-\alpha+2}}dy.\nonumber
\end{align}
By the symmetry of the double integral,
\begin{align}
2\int_\Omega dx\int_\Omega\frac{(x-y)\cdot x u(x)u(y)}{|x-y|^{N-\alpha+2}}dy &= \int_\Omega dx\int_\Omega\frac{(x-y)\cdot x u(x)u(y)}{|x-y|^{N-\alpha+2}}dy\nonumber\\
&\qquad - \int_\Omega dx\int_\Omega\frac{(x-y)\cdot y u(x)u(y)}{|x-y|^{N-\alpha+2}}dy\nonumber\\
&=\int_\Omega dx\int_\Omega\frac{u(x)u(y)}{|x-y|^{N-\alpha}}dy.\nonumber
\end{align}
Hence,
\begin{align}
\int_\Omega(I_\alpha*_\Omega u)(x\cdot\nabla u)dx =-\frac{N+\alpha}{2} C_{N,\alpha} \int_\Omega dx\int_\Omega\frac{u(x)u(y)}{|x-y|^{N-\alpha}}dy.\nonumber
\end{align}

\begin{flushright}
$\Box$
\end{flushright}

Now we prove Pohozaev identity.
\begin{lemma}\label{l:Pohozaev}
For any solution to Problem (\ref{e:000}), it holds that
\begin{align}
0&=\int_{\partial\Omega}\Big|\frac{\partial u}{\partial n}\Big|^2x\cdot d\vec{S}_x + (N-2)\int_\Omega|\nabla u|^2dx- \lambda C_{N,\alpha}(N+\alpha) \int_\Omega dy\int_\Omega\frac{u(x)u(y)}{|x-y|^{N-\alpha}}dx\nonumber\\
&\quad  -(N-2)\int_\Omega|u|^{2^*}dx.\nonumber
\end{align}
\end{lemma}
\noindent{\bf Proof.}
To do this, only need to notice that
\begin{align}
0=\int_\Omega(\Delta u+\lambda I_\alpha*_\Omega u+|u|^{2^*-2}u)(x\cdot\nabla u)=I_1+I_2+I_3.\nonumber
\end{align}
Here,
\begin{align}
I_1&=\int_\Omega\Delta u(x\cdot\nabla u)dx\nonumber\\
&=\int_\Omega\Big[\mbox{div}\big(\nabla u(x\cdot\nabla u)\big)-|\nabla u|^2-x\cdot\nabla\big( \frac{|\nabla u|^2}{2}\big)\Big]dx\nonumber\\
&=\int_\Omega \Big[ \mbox{div}\big(\nabla(x\cdot\nabla u)-x\frac{|\nabla u|^2}{2}\big) +\frac{N-2}{2}|\nabla u|^2\Big]dx.\nonumber
\end{align}
Since for any $x\in\Omega$, $x\cdot\nabla u=x\cdot n\frac{\partial u}{\partial n}$, then
\begin{align}
I_1=\frac{1}{2}\int_{\partial\Omega}\Big|\frac{\partial u}{\partial n}\Big|^2x\cdot d\vec{S}_x+ \frac{N-2}{2}\int_\Omega|\nabla u|^2dx.\nonumber
\end{align}
Moreover,
\begin{align}
I_2=\lambda\int_\Omega(I_\alpha*_\Omega u)(x\cdot\nabla u)dx=-\frac{\lambda C_{N,\alpha}(N+\alpha)}{2} \int_\Omega dx\int_\Omega\frac{u(x)u(y)}{|x-y|^{N-\alpha}}dy.\nonumber
\end{align}
due to Lemma \ref{l:PotentialPohozaev}. And,
\begin{align}
I_3&=\int_\Omega|u|^{2^*-2}u\cdot(x\cdot\nabla u)dx=\int_\Omega\Big[\mbox{div}\big(\frac{x}{2^*} |u|^{2^*}\big)-\frac{N-2}{2}|u|^{2^*}\Big]dx\nonumber\\
&=-\frac{N-2}{2}\int_\Omega|u|^{2^*}dx.\nonumber
\end{align}
The lemmas follows immediately.

\begin{flushright}
$\Box$
\end{flushright}

\subsection{Proof of Theorem \ref{t:nonexistence}}
In this part, we prove Theorem \ref{t:nonexistence}.

\vs
\noindent{\bf Proof of Theorem \ref{t:nonexistence}.}
Let $u$ be a solution to Problem (\ref{e:000}), then
\begin{align}
\int_\Omega|\nabla u|^2 dx-\lambda C_{N,\alpha}\int_\Omega dx \int_\Omega\frac{u(x)u(y)}{|x-y|^{N-\alpha}}dy =\int_\Omega|u|^{2^*}dx.\nonumber
\end{align}
By Lemma \ref{l:Pohozaev}, we get
\begin{align}
0= \int_{\partial\Omega}\Big|\frac{\partial u}{\partial n}\Big|^2x\cdot d\vec{S}_x-(\alpha+2)\lambda C_{N,\alpha} \int_\Omega dx \int_\Omega\frac{u(x)u(y)}{|x-y|^{N-\alpha}}dy.
\end{align}
Suppose that $\Omega$ is star-shaped, if $\lambda<0$, we get $u\equiv0$ immediately. If $\lambda=0$, $u\equiv 0$ due to the unique continuation, cf. \cite{Heinz1995}.

\begin{flushright}
$\Box$
\end{flushright}

\subsection{Proof of Theorem \ref{t:nonexistence2}}
In this part, we prove Theorem \ref{t:nonexistence2}.

\vs
\noindent{\bf Proof of Theorem \ref{t:nonexistence2}.}
By Pohozaev identity, if Problem (\ref{e:000}) admits a positive solution,
\begin{align}
\lambda C_{N,\alpha}\int_{\Omega}dx\int_\Omega\frac{u(x)u(y)}{|x-y|^{N-\alpha}}dy & =\frac{1}{\alpha+2}\int_{\partial\Omega}\Big|\frac{\partial u}{\partial n}\Big|^2x\cdot d\vec{S}_x\nonumber\\
&\geq \frac{c_0}{\alpha+2} \int_{\partial\Omega}\Big|\frac{\partial u}{\partial n}\Big|^2 dS_x\nonumber\\
&\geq\frac{c_0|\partial\Omega|_{N-1}^{-1}}{\alpha+2}\Big(\int_{\partial\Omega}\frac{\partial u}{\partial n} dS_x\Big)^2\nonumber\\
&=\frac{c_0|\partial\Omega|_{N-1}^{-1}}{\alpha+2}\Big(\int_\Omega(-\Delta u)dx\Big)^2.\nonumber
\end{align}
Noticing that $u|_{\partial\Omega}=0$, by the Sobolev inequality involving $|\Delta u|_1$ (see, for instance, \cite[Theorem 8]{BrezisStrauss1973}), if $N-\alpha-4<0$, there exists a constant $C_*>0$ independence in $u$ and
\begin{align}
|u|_{\frac{2N}{N+\alpha}}\leq C\|u\|_{W^{1,q}(\Omega)}\leq C_* |\Delta u|_1\nonumber
\end{align}
for any $1\leq q<\frac{N}{N-1}$.
Then, by Hardy-Littlewood-Sobolev inequality, there exists a constant $C=C(N,\alpha)>0$ such that
\begin{align}
\lambda C_{N,\alpha}\int_{\Omega}dx\int_\Omega\frac{u(x)u(y)}{|x-y|^{N-\alpha}}dy \geq \frac{C|\Omega|_N^{\frac{N-\alpha-4}{N}}}{|\partial\Omega|_{N-1}} C_{N,\alpha} \int_{\Omega}dx\int_\Omega\frac{u(x)u(y)}{|x-y|^{N-\alpha}}dy.\nonumber
\end{align}
This proves the result.

\begin{flushright}
$\Box$
\end{flushright}

\section{Proof of Theorems \ref{t:existence} and \ref{t:multiplicity}: Existence and multiplicity results}\label{Sec:Existence}
\subsection{Estimate of the mountain-pass level}
The proofs of Theorems \ref{t:existence} and \ref{t:multiplicity} are based on \cite{BrezisNirenberg1983,CeramiFortStruwe1984}.
To begin with,
we use the Talenti instanton to test the mountain pass level of functional $I_\alpha$. Talenti instanton is defined as
\begin{align}
U_{x_0,\mu}(x)=C_0\frac{\mu^{\frac{N-2}{2}}}{(1+\mu^2|x-x_0|^2)^\frac{N-2}{2}}\nonumber
\end{align}
with $C_0$ is a positive constant depends on $N$ only, cf. \cite[pp. 18]{CaoPengYan2021}.
A minimizing formulation of the mountain pass level is that
\begin{align}
m_{\alpha,\lambda}=\inf\Big\{J_{\alpha,\lambda}(u)\Big|\int_\Omega|u|^{2^*}dx=1\Big\}\nonumber
\end{align}
with
\begin{align}
J_{\alpha,\lambda}(u)=\int_\Omega|\nabla u|^2 dx-\lambda\mathcal{D}_\alpha(u).\nonumber
\end{align}
Consider the minimization
\begin{align}
S=\inf\Big\{\int_{\mathbb{R}^N}|\nabla u|^2dx\Big|\int_{\mathbb{R}^N}|u|^{2^*}dx=1\Big\}.\nonumber
\end{align}
Without loss of generality, let us assume that $\int_{\mathbb{R}^N}|U_{x_0,\mu}|^{2^*}dx=1$ and that $U_{x_0,\mu}$ attains $S$.
Moreover, let us assume that $0\in\Omega$. For a small positive number $\delta>0$, define a smooth cutoff function
\begin{equation}
\xi(x)=\left\{
\begin{aligned}
1, & |x|<\delta, \nonumber\\
0, & |x|>2\delta.\nonumber
\end{aligned}
\right.
\end{equation}
By a routine computation as in \cite{BrezisNirenberg1983} and in \cite[pp. 19]{CaoPengYan2021},
\begin{proposition}
It holds that
\begin{itemize}
  \item [$(1).$] $\int_\Omega|\xi U_{0,\mu}|^{2^*}=\int_{\mathbb{R}^N}U_{0,1}^{2^*}+O(\frac{1}{\mu^N})$;
  \item [$(2).$] $\int_\Omega|\nabla(\xi U_{0,\mu})|^2=\int_{\mathbb{R}^N}|\nabla U_{0,1}|^2+ O(\frac{1}{\mu^{N-2}})$.
\end{itemize}
\end{proposition}
Now we will give an estimate on $C_{N,\alpha}\int_\Omega dx\int_\Omega\frac{\xi(x)\xi(y) U_{0,\mu}(x)U_{0,\mu}(y)}{|x-y|^{N-\alpha}}dy$. To this end, we first study its limit, i.e. the integral
\begin{align}\label{inte:limit}
C_{N,\alpha}\int_{\mathbb{R}^N} dx\int_{\mathbb{R}^N}\frac{U_{0,1}(x)U_{0,1}(y)}{|x-y|^{N-\alpha}}dy.
\end{align}
\begin{lemma}
Integral (\ref{inte:limit}) converges if and only if $N-\alpha-4>0$.
\end{lemma}
\noindent{\bf Proof.}
Notice that
\begin{align}
\mbox{(\ref{inte:limit})}=|I_\frac{\alpha}{2}*U_{0,1}|_2^2=C||\xi|^{-\frac{\alpha}{2}} \widehat{U_{0,1}}(\xi)|_2^2. \nonumber
\end{align}
Now we consider $\widehat{U_{0,1}}(\xi)$, the Fourier transform of $U_{0,1}$. Do as in \cite[pp. 360]{Lieb1983},
\begin{align}
\widehat{\Big[\Big(\frac{1}{1+|x|^2}\Big)^\frac{N-2}{2}\Big]}(\xi)= \frac{\pi^\frac{N}{2}2^{\frac{2}{N} -\frac{N}{2}+2}}{\Gamma\Big(\frac{N}{2}-1\Big)} |\xi|^{-1}K_1(|\xi|).\nonumber
\end{align}
Here, $K_1$ is the modified Bessel function of the second kind $K_\upsilon$ with $\upsilon=1$.
Using \cite[(14)/pp. 88]{Watson1995}, we get
\begin{align}\label{Asy:0}
|\xi|^{-\alpha}|\widehat{U_{0,1}}(\xi)|^2|\xi|^{N-1}\sim C|\xi|^{(N-\alpha-4)-1}
\end{align}
as $|\xi|\to 0$.
By \cite[(1)/pp. 202]{Watson1995}, we get
\begin{align}\label{Asy:infty}
|\xi|^{-\alpha}|\widehat{U_{0,1}}(\xi)|^2|\xi|^{N-1}\sim C e^{-|\xi|}|\xi|^{N-\alpha-4}
\end{align}
as $|\xi|\to\infty$.
Combining (\ref{Asy:0}) and (\ref{Asy:infty}), we prove the lemma.
\begin{flushright}
$\Box$
\end{flushright}

We only need to check
$\mathcal{D}_\alpha(\xi U_{0,\mu})$.
By an easy computation, we get
\begin{align}
\mathcal{D}_\alpha(\xi U_{0,\mu}) & =C_{N,\alpha}\int_\Omega dy\int_\Omega \frac{\xi(x)\xi(y)U_{0,\mu}(x)U_{0,\mu}(y)}{|x-y|^{N-\alpha}}\nonumber\\
&=C_{N,\alpha}\mu^{-2-\alpha}\int_{\mathbb{R}^N\times\mathbb{R}^N} \frac{\xi\big(\frac{x}{\mu}\big)\xi\big(\frac{x}{\mu}\big)U_{0,1}(x)U_{0,1}(y)}{|x-y|^{N-\alpha}}dxdy. \nonumber
\end{align}
Denote
\begin{align}
D_R:=C_{N,\alpha}\int_{R_R(0)\times B_R(0)} \frac{U_{0,1}(x)U_{0,1}(y)}{|x-y|^{N-\alpha}}dxdy. \nonumber
\end{align}
\begin{lemma}\label{l:DR}
It holds that
\begin{itemize}
  \item [$(1).$] If $N-\alpha-4>0$, $D_R\to C$ for some $C>0$ as $R\to\infty$;
  \item [$(2).$] If $N-\alpha-4\leq0$, $D_R\to+\infty$ as $R\to\infty$;
  \item [$(3).$] Given $R>0$, for large $\mu>0$, $\mathcal{D}_\alpha(\xi U_{0,\mu})\geq D_R\mu^{-2-\alpha}$.
\end{itemize}
\end{lemma}
The next lemma follows immediately.
\begin{lemma}\label{l:EnergyLevel}
$m_{\alpha,\lambda}<S$ if one of the following holds:
\begin{itemize}
  \item [$(1).$] $N-\alpha-4\geq0$ and any $\lambda>0$;
  \item [$(2).$] $N-\alpha-4<0$ and $\lambda>0$ large.
\end{itemize}
\end{lemma}
\noindent{\bf Proof.}
Denote $v_\mu=\frac{\xi U_{0,\mu}}{|\xi U_{0,\mu}|_{2^*}}$. Then,
\begin{align}
m_{\alpha,\lambda}&\leq I_\lambda(v_\mu)= \frac{\int_\Omega|\nabla(\xi U_{0,\mu})|^2-\lambda\mathcal{D}_\alpha(\xi U_{0,\mu})}{(\int_\Omega|\xi U_{0,\mu}|^{2^*})^{\frac{N-2}{N}}}\nonumber\\
&\leq\frac{\int_{\mathbb{R}^N}|\nabla U_{0,1}|^2+O(\mu^{2-N})-\lambda D_R \mu^{-2-\alpha}}{(1+O(\mu^{-N}))^{\frac{N-2}{N}}}=\int_{\mathbb{R}^N}|\nabla U_{0,1}|^2+ \frac{O(\mu^{2-N})- \lambda D_R \mu^{-2-\alpha}}{1+O(\mu^{2-N})}\nonumber\\
&=S+T_{\lambda,\mu,\alpha,N}.\nonumber
\end{align}
Here,
\begin{align}
T_{\lambda,\mu,\alpha,N}=\frac{O(\mu^{2-N})- \lambda D_R \mu^{-2-\alpha}}{1+O(\mu^{2-N})}.\nonumber
\end{align}

In order to prove $m_{\alpha,\lambda}\leq I_\lambda(v_\mu)<S$, it is sufficient to check that $T_{\lambda,\mu,\alpha,N}<0$.

\vs

\noindent{\bf Case. 1. }$N-\alpha-4>0$.

\vs

If $N-\alpha-4>0$, for any $\lambda>0$, letting $\mu>0$ be sufficiently large, we get $T_{\lambda,\mu,\alpha,N}<0$.

\vs

\noindent{\bf Case. 2. }$N-\alpha-4=0$.

\vs

If $N-\alpha-4=0$, for any $\lambda>0$, Lemma \ref{l:DR} implies that $D_R\to+\infty$ as $R\to\infty$. Then, for large $R>0$ and $\mu>0$, we get $T_{\lambda,\mu,\alpha,N}<0$ again.

\vs

\noindent{\bf Case. 3. }$N-\alpha-4<0$.

\vs

If $N-\alpha-4<0$, for arbitrarily fixed $\mu$ and $R$, let $\lambda>0$ be sufficiently large, we get $T_{\lambda,\mu,\alpha,N}<0$.

\begin{flushright}
$\Box$
\end{flushright}

\subsection{A local compactness}
In this subsection, we prove a convergence theorem for the $(PS)$ sequence of the energy functional
\begin{align}
E_{\alpha,\lambda}(u)=\frac{1}{2}\int_\Omega|\nabla u|^2 -\frac{\lambda}{2}\mathcal{D}_\alpha(u) -\frac{1}{2^*}\int_\Omega|u|^{2^*}\nonumber
\end{align}
satisfying certain assumptions. To be precise, we get
\begin{lemma}\label{l:compactness}
For any sequence $\{u_n\}_n\in H_0^1(\Omega)$ such that
\begin{itemize}
  \item [$(1).$] $E_{\alpha,\lambda}(u_n)\to c<\frac{1}{N}S^\frac{N}{2}$;
  \item [$(2).$] $\nabla E_{\alpha,\lambda}(u_n)\to0$ in $H^{-1}(\Omega)$,
\end{itemize}
then $\{u_k\}_k$ admits a convergent subsequence in $H_0^1(\Omega)$.
\end{lemma}
\noindent{\bf Proof.}
We prove it in steps.

\vs

\noindent{\bf Step 1. }$\{u_n\}_n$ is bounded in $H^1_0(\Omega)$.

\vs

It is evident that
\begin{align}
o_n(1)\|u_n\|+c+1\geq\frac{1}{N}|u_n|_{2^*}^{2^*}.\nonumber
\end{align}
Then,
\begin{align}
|\nabla u_n|_2^2 & \leq C+o_n(1)\|u_n\|+\frac{|\lambda|}{N}\mathcal{D}_\alpha(u_n)\nonumber\\
&\leq C+o_n(1)\|u_n\|+(C+o_n(1)\|u_n\|)^\frac{N-2}{N}.\nonumber
\end{align}
This implies the boundedness of $\{u_n\}_n$ in $H^1_0(\Omega)$.
\vs

\noindent{\bf Step 2. }$\{u_n\}_n$ converges in $H^1_0(\Omega)$.

\vs

First notice that there exists a $u_\infty\in H_0^1(\Omega)$ such that $u_n\rightharpoonup u_\infty$ in $H_0^1(\Omega)$. We will prove that $u_\infty\to u_\infty$ in $H_0^1(\Omega)$. Then, it is evident that
\begin{itemize}
  \item $\nabla E_{\alpha,\lambda}(u_\infty)=0$;
  \item $\int_\Omega|\nabla u_n|^2 dx=\int_\Omega|\nabla(u_n- u_\infty)|^2 dx+\int_\Omega|\nabla u_\infty|^2 dx+o_n(1)$;
  \item $\mathcal{D}_\alpha(u_n)\to\mathcal{D}_\alpha(u_\infty)$;
  \item $\int_\Omega|u_n|^{2^*} dx=\int_\Omega|u_n -u_\infty|^{2^*} dx+\int_\Omega|u_\infty|^{2^*} dx +o_n(1)$.
\end{itemize}
Here, the last assertion applied Brezis-Lieb theorem. On the other hand, notice that
\begin{align}
o_n(1) & =\nabla E_{\alpha,\lambda}(u_n)(u_n-u_\infty)\nonumber\\
&=\Big(\nabla E_{\alpha,\lambda}(u_n)-\nabla E_{\alpha,\lambda}(u_\infty)\Big)(u_n-u_\infty)\nonumber\\
&=\int_\Omega|\nabla(u_n-u_\infty)|^2 dx-\int_\Omega|u_n|^{2^*-2}u_n(u_n -u_\infty) dx+o_n(1).\nonumber
\end{align}
Here, using Brezis-Lieb theorem,
\begin{align}
\int_\Omega|u_n|^{2^*-2}u_n(u_n -u_\infty) dx =\int_\Omega|u_n|^{2^*} dx-\int_\Omega|u_\infty|^{2^*} dx= \int_\Omega|u_n -u_\infty|^{2^*}dx+o_n(1).\nonumber
\end{align}
Then, we get
\begin{align}
\int_\Omega|\nabla (u_n-u_\infty)|^2 dx=\int_\Omega|u_n-u_\infty|^{2^*} dx+o_n(1).\nonumber
\end{align}
It follows that
\begin{align}
E_{\alpha,0}(u_n-u_\infty)=\frac{1}{N}\int_\Omega|\nabla (u_n-u_\infty)|^2 dx +o_n(1).\nonumber
\end{align}
Then,
\begin{align}\label{equa:convergence1}
\int_\Omega|\nabla(u_n-u_\infty)|^2 dx =N(E_{\alpha,\lambda}(u_n)-E_{\alpha,\lambda}(u_\infty))+o_n(1) <S^\frac{N}{2}.
\end{align}
On the other hand, notice that
\begin{align}\label{equa:convergence2}
\|u_n -u_\infty\|^2(1-S^{-\frac{2^*}{2}}\|u_n -u_\infty\|^{2^*-2})\leq\int_\Omega|\nabla (u_n-u_\infty)|^2 dx  -\int_\Omega|u_n -u_\infty|^{2^*} dx=o_n(1).
\end{align}
Combining (\ref{equa:convergence1}) and (\ref{equa:convergence2}), we get $\|u_n-u_\infty\|\to0$.
\begin{flushright}
$\Box$
\end{flushright}

\subsection{A maximum principle for a nonlocal operator}
In this part, we prove the following maximum principle for the nonlocal operator
\begin{align}
S_\lambda(u):=-\Delta u-\lambda I_{\alpha}*_\Omega u.\nonumber
\end{align}
\begin{lemma}\label{l:maximumprinciple}
Suppose $0<\lambda<\lambda_1$ and $u\in H_0^1(\Omega)\backslash \{0\}$ satisfies that $S_\lambda(u)\geq0$, then $u\geq0$ in $\Omega$. Here,
\begin{align}
\lambda_1 =\inf_{u\neq 0}\frac{\int_\Omega|\nabla u|^2 dx}{\mathcal{D}_\alpha(u)}.\nonumber
\end{align}
\end{lemma}
\noindent{\bf Proof.}
According to the assumption, it holds that
\begin{align}\label{ineq:Pointwise}
-\Delta u-\lambda I_\alpha*_\Omega u\geq 0
\end{align}
Testing (\ref{ineq:Pointwise}) with $u_-(x)=\min\{0,u(x)\}$, we get
\begin{align}\label{ineq:gradient}
0&\geq\int_\Omega|\nabla u|^2 dx -\lambda C_{N,\alpha}\int_\Omega\int_\Omega\frac{u_+(x)u_-(y)}{|x-y|^{N-\alpha}}dxdy -\lambda\int_\Omega \int_\Omega\frac{u_-(x)u_-(y)}{|x-y|^{N-\alpha}}dxdy \nonumber\\
&\geq\int_\Omega|\nabla u_-|^2 dx -\lambda C_{N,\alpha}\int_\Omega \int_\Omega\frac{u_-(x)u_-(y)}{|x-y|^{N-\alpha}}dxdy.
\end{align}
with $u_+(x)=\max\{0,u(x)\}$. Now we study $\lambda C_{N,\alpha}\int_\Omega \int_\Omega\frac{u_-(x)u_-(y)}{|x-y|^{N-\alpha}}dxdy$.
By Hardy-Littlewood-Sobolev inequality and Sobolev inequality, we get
\begin{align}
\lambda C_{N,\alpha}\int_\Omega \int_\Omega\frac{u_-(x)u_-(y)}{|x-y|^{N-\alpha}}dxdy &=\lambda\mathcal{D}_\alpha(u)\leq\lambda\cdot\lambda_1^{-1}\int_\Omega|\nabla u_-|^2 dx.\nonumber
\end{align}
Here,
\begin{align}
\lambda_1=\inf_{u\in\mathcal{M}}\int_\Omega|\nabla u|^2 =\inf_{u\neq 0}\frac{\int|\nabla u|^2 dx}{\mathcal{D}_\alpha(u)}.\nonumber
\end{align}
Plugging this into (\ref{ineq:gradient}), we get
\begin{align}
0\geq|\nabla u_-|_2^2-\lambda\cdot \lambda_1^{-1}|\nabla u_-|_2^2=\Big(1-\frac{\lambda}{\lambda_1}\Big) |\nabla u_-|_2^2.\nonumber
\end{align}
Notice that $0<\lambda<\lambda_1$. We get $u_-=0$ and complete the proof.
\begin{flushright}
$\Box$
\end{flushright}

\subsection{Proof of Theorem \ref{t:existence}}
In this part, we prove Theorem \ref{t:existence}.

\vs

\noindent{\bf Proof of Theorem \ref{t:existence}.}
By Lemma \ref{l:EnergyLevel}, for the cases
\begin{itemize}
  \item [$(1).$] $N-\alpha-4\geq0$ and any $\lambda>0$;
  \item [$(2).$] $N-\alpha-4<0$ and $\lambda>0$ large,
\end{itemize}
there exists a a sequence $\{u_n\}_n\subset H_0^1(\Omega)$ such that as $n\to\infty$,
\begin{align}
J_{\alpha,\lambda}(u_n)\to m_{\alpha,\lambda}=\inf\Big\{\int_\Omega|\nabla u|^2 dx-\lambda\mathcal{D}_{\alpha}(u)\Big|\int_\Omega|u|^{2^*}dx=1\Big\} <S.
\end{align}
First, we show that $\{u_n\}_n$ is bounded in $H_0^1(\Omega)$. Indeed, by a direct computation,
\begin{align}
m_{\alpha,\lambda}\geq\int_\Omega|\nabla u_n|^2 dx-\lambda\mathcal{D}_\alpha(u_n)\geq\int_\Omega|\nabla u_n|^2 dx -C|u_n|_{2^*}^2\geq \int_\Omega|\nabla u_n|^2 dx -C.\nonumber
\end{align}
Then, $\{u_n\}_n$ is bounded in $H_0^1(\Omega)$ and therefore there exists a $u\in H_0^1(\Omega)$ such that $u_n\rightharpoonup u$ in $H_0^1(\Omega)$.

On the other hand, by Brezis-Lieb theorem and the weak convergence of $\{u_n\}_n$ to $u$, we get
\begin{align}
\int_\Omega|u_n|^{2^*}dx-\int_\Omega|u_n -u|^{2^*}dx=\int_\Omega|u|^{2^*}dx+o_n(1),\nonumber
\end{align}
\begin{align}
\int_\Omega|\nabla u_n|^{2}dx-\int_\Omega|\nabla (u_n -u)|^{2}dx=\int_\Omega|\nabla u|^{2}dx+o_n(1)\nonumber
\end{align}
and
\begin{align}
\mathcal{D}_\alpha(u_m)=\mathcal{D}_\alpha(u)+o_n(1).\nonumber
\end{align}
Then,
\begin{align}
m_{\alpha,\lambda}+o_n(1)&=\int_\Omega|\nabla(u_n -u)|^2 dx+\int_\Omega|\nabla u|^2 dx-\lambda\mathcal{D}_\alpha(u)+o_n(1)\nonumber\\
&\geq S|u_n -u|_{2^*}^2+m_{\alpha,\lambda}|u|_{2^*}^2+o_n(1)\nonumber\\
&\geq S|u_n -u|_{2^*}^{2^*}+m_{\alpha,\lambda}|u|_{2^*}^{2^*}+o_n(1)\nonumber\\
&\geq(S-m_{\alpha,\lambda})|u_n -u|_{2^*}^{2^*}+m_{\alpha,\lambda}+o_n(1).\nonumber
\end{align}
Since $m_{\alpha,\lambda}<S$, $u_n\to u$ in $L^{2^*}(\Omega)$ and $|u|_{2^*}=1$. Therefore,
\begin{align}
m_{\alpha,\lambda}\leq J_{\alpha,\lambda}(u)\leq\lim_{n\to\infty}J_{\alpha,\lambda}(u_n)=m_{\alpha,\lambda},\nonumber
\end{align}
which proves that $u$ attains the minimum $m_{\alpha,\lambda}$. This implies that $u$ is a solution to Problem (\ref{e:000}).

To prove the second part, it is necessary to modify the settings as follows.
Denote
\begin{align}
m^+_{\alpha,\lambda}=\inf\Big\{J_{\alpha,\lambda}(u)\Big|\int_\Omega|u_+|^{2^*}dx=1\Big\}.\nonumber
\end{align}
Notice that the set $\Big\{u\in H_0^1(\Omega)\Big|\int_\Omega |u_+|^{2^*}dx =1\Big\}$ is a $C^1$ manifold. Then the above computation can be repeated. To be precise, we get
\begin{claim}
$m_{\alpha,\lambda}^+<S$ if one of the following holds:
\begin{itemize}
  \item [$(1).$] $N-\alpha-4\geq0$ and any $\lambda>0$;
  \item [$(2).$] $N-\alpha-4<0$ and $\lambda>0$ large.
\end{itemize}
\end{claim}
Then, there exists a minimizing sequence $\{u'_n\}_n$ of $m_{\alpha,\lambda}^+$.
Repeating the argument as above, we find a non-negative minimizer $u'$ of $m_{\alpha,\lambda}^+$, which solves the problem
\begin{equation}
    \left\{
   \begin{array}{lr}
     -{\Delta}u=\lambda I_\alpha*_\Omega u+|u_+|^{2^*-2}u_+\mbox{ in }\Omega ,\nonumber\\
     u\in H_0^1(\Omega).\nonumber
   \end{array}
   \right.
\end{equation}
Since $0<\lambda<\lambda_1$, by Lemma \ref{l:maximumprinciple} and Lemma \ref{l:RegularNonlinearC2a}, $u'$ is a positive solution to Problem (\ref{e:000}).

\begin{flushright}
$\Box$
\end{flushright}

\subsection{Proof of Theorem \ref{t:multiplicity}}\label{Subsec:multiple}
We follow the notation in the last subsection.

\vs

\noindent{\bf Proof of Theorem \ref{t:multiplicity}.}
We apply the idea in \cite{CeramiFortStruwe1984}.
Let $A_j=\mbox{span}\{\phi_1,\cdots,\phi_j\}\cap\{u\in H_0^1(\Omega)||u|_{2^*}=1\}$. By a direct computation,
\begin{align}
\sup_{u\in A_j}\frac{\int_\Omega|\nabla u|^2 dx -\lambda\mathcal{D}_\alpha(u)}{|u|_{2^*}^2} \leq (\lambda_j-\lambda)\sup_{A_j}\frac{\mathcal{D}_\alpha(u)}{|u|_{2^*}^2} <\upsilon\cdot C_{HLS}C_{N,\alpha}|\Omega|_N^\frac{2+\alpha}{N}.\nonumber
\end{align}
Here, $C_{HLS}$ is the best constant in the Hardy-Littlewood-Sobolev inequality
\begin{align}
\int_{\mathbb{R}^N} dy\int_{\mathbb{R}^N}\frac{u(x)v(y)}{|x-y|^{N-\alpha}}dx \leq C_{HLS}|u|_\frac{2N}{N+\alpha}|v|_\frac{2N}{N+\alpha}.\nonumber
\end{align}
By the assumption in Theorem \ref{t:multiplicity},
\begin{align}
\sup_{u\in A_j}\frac{\int_\Omega|\nabla u|^2 dx -\lambda\mathcal{D}_\alpha(u)}{|u|_{2^*}^2} \leq (\lambda_j-\lambda)\sup_{A_j}\frac{\mathcal{D}_\alpha(u)}{|u|_{2^*}^2} <S.\nonumber
\end{align}
Denote
\begin{align}
F_{\alpha,\lambda}(u)=\frac{\int_\Omega|\nabla u|^2 dx -\lambda\mathcal{D}_\alpha(u)}{|u|_{2^*}^2}.\nonumber
\end{align}
Then, $\gamma(\mbox{int}(F_{\alpha,\lambda}^{S})\cap \{u\in H_0^1(\Omega)||u|_{2^*}=1\})\geq j+m$.
Denote $\mathcal{L}=\{u\in H_0^1(\Omega)|\int_\Omega|u|^{2^*}dx=1\}$, $L:=\overline{\mbox{span}\{\phi_j,\phi_{j+1},\cdots\}}\cap\mathcal{L}$ and $\Gamma_k=\{A\subset\mathcal{L}|\gamma(A)\geq j+k\}$ for $k=1,\cdots,m$.
Here, $\gamma$ is the $\mathbb{Z}_2$ genus generated by the anti-podal action $u\mapsto-u$.
Define the minimax values
\begin{align}
c_k=\inf_{A\in\Gamma_k}\sup_{x\in A}J_{\alpha,\lambda}(u).
\end{align}
Notice that for any $k=1,\cdots,m$ and any $A\in\Gamma_k$, $A\cap L\neq\emptyset$. By the intersection lemma \cite[Proposition 7.8]{RabinowitzBook1986}, we get
\begin{align}
c_k=\inf_{A\in\Gamma_k}\sup_{u\in A}J_{\alpha,\lambda}(u)\geq\inf_{u\in L}J_{\alpha,\lambda}(u)\geq(\lambda_j-\lambda)\int_\Omega|\nabla u|^2 dx\Big|_{\mathcal{L}}\geq c_*>0.\nonumber
\end{align}
By \cite{AmbrosettiRabinowitz1973} and \cite[Theorem 8.10]{RabinowitzBook1986}, $c_k$ are critical values of $J_{\alpha,\lambda}$ for $k=1,\cdots,m$.
For any $k=1,\cdots,m$, there exists a $(PS)$ sequence $\{u_n^k\}_n$ of $J_{\alpha,\lambda}(u)|_{\mathcal{L}}$ such that $J_{\alpha,\lambda}(u_n^k)\to c_k$ as $n\to\infty$. Then, the sequence $\{\widetilde{u}_n^k\}_n$ is a $(PS)$ sequence of $E_{\alpha,\lambda}$ such that $\lim_{n\to\infty}E_{\alpha,\lambda}(\widetilde{u}_n^k)<\frac{1}{N} S^\frac{N}{2}$. By Lemma \ref{l:compactness}, $\{\widetilde{u}_n^k\}_n$ converges to a critical point of $E_{\alpha,\lambda}$ and Theorem \ref{t:multiplicity} follows immediately.

\begin{flushright}
$\Box$
\end{flushright}

\Vs

\Vs

\noindent{\bf Acknowledgment. }
Both of the authors would like to express their sincere appreciation to the anonymous reviewers for their helpful suggestions on improving this paper.

\Vs

\noindent{\bf Fundings. }
In this research, HL was supported by FAPESP PROC 2022/15812-0.

\Vs

\noindent{\bf Competing Interests. }
The authors have no relevant financial or non-financial competing interests.

\Vs

\noindent{\bf Data availability statement. }
No datasets were generated or analysed during the current study.

\Vs\Vs


{\footnotesize

\begin {thebibliography}{44}

\bibitem{AlvesCassaniTarsiYang2016}
Alves, C. O, Cassani, D, Tarsi, C, Yang, M,
Existence and concentration of ground state solutions for a critical nonlocal Schr\"odinger equation in $\mathbb{R}^2$.
J. Differ. Equations 261, No. 3, 1933-1972 (2016).



\bibitem{AmbrosettiRabinowitz1973}
Ambrosetti, A, Rabinowitz, P. H,
Dual variational methods in critical point theory and applications.
J. Funct. Anal. 14, 349-381 (1973).





\bibitem{BahGroDonBas2014}
Bahrami, M, Gro$\beta$ardt, A, Donadi, S, Bassi, A,
The Schr\"odinger-Newton equation and its foundations.
New J. Phys. 16, No. 11, Article ID 115007, 17 p. (2014).




\bibitem{BellazziniFrankVisciglia2014}
Bellazzini, J, Frank, R. L, Visciglia, N,
Maximizers for Gagliardo-Nirenberg inequalities and related non-local problems.
Math. Ann. 360, No. 3-4, 653-673 (2014).



\bibitem{BenciFort1998}
Benci, V, Fortunato, D,
An eigenvalue problem for the Schr\"odinger-Maxwell equations.
Topol. Methods Nonlinear Anal. 11, No. 2, 283-293 (1998).



\bibitem{BokaLopezSoler2003}
Bokanowski, O, L\'opez, Jos\'e L, Soler, J,
On an exchange interaction model for quantum transport: the Schr\"odinger-Poisson-Slater system.
Math. Models Methods Appl. Sci. 13, No. 10, 1397-1412 (2003).


\bibitem{BrezisNirenberg1983}
Br\'ezis, H, Nirenberg, L,
Positive solutions of nonlinear elliptic equations involving critical Sobolev exponents.
Commun. Pure Appl. Math. 36, 437-477 (1983).


\bibitem{BrezisStrauss1973}
Br\'ezis, H, Strauss, W. A,
Semi-linear second-order elliptic equations in $L^1$.
J. Math. Soc. Japan 25, 565-590 (1973).


\bibitem{CaoPengYan2021}
Cao, D, Peng, S, Yan, S,
Singularly perturbed methods for nonlinear elliptic problems.
Cambridge Studies in Advanced Mathematics 191. Cambridge: Cambridge University Press. ix, 252 p. (2021).

\bibitem{CeramiFortStruwe1984}
Cerami, G, Fortunato, D, Struwe, M,
Bifurcation and multiplicity results for nonlinear elliptic problems involving critical Sobolev exponents.
Ann. Inst. Henri Poincar\'e, Anal. Non Lin\'eaire 1, 341-350 (1984).



\bibitem{CeramiMolle2019}
Cerami, G, Molle, R,
Multiple positive bound states for critical Schr\"odinger-Poisson systems.
ESAIM, Control Optim. Calc. Var. 25, Paper No. 73, 29 p. (2019).




\bibitem{ChenJinLiLim2005}
Chen, W, Jin, C, Li, C, Lim, J,
Weighted Hardy-Littlewood-Sobolev inequalities and systems of integral equations. (English) Zbl 1147.45301
Discrete Contin. Dyn. Syst. 2005, Suppl., 164-172 (2005).




\bibitem{ChenLi2008}
Chen, W, Li, C,
The best constant in a weighted Hardy-Littlewood-Sobolev inequality.
Proc. Am. Math. Soc. 136, No. 3, 955-962 (2008).




\bibitem{ChenLiOu2005}
Chen, W, Li, C, Ou, B,
Classification of solutions for a system of integral equations.
Commun. Partial Differ. Equations 30, No. 1, 59-65 (2005).




\bibitem{DouZhu20151}
Dou, J, Zhu, M,
Reversed Hardy-Littewood-Sobolev inequality.
Int. Math. Res. Not. 2015, No. 19, 9696-9726 (2015).



\bibitem{DouZhu20152}
Dou, J, Zhu, M,
Sharp Hardy-Littlewood-Sobolev inequality on the upper half space.
Int. Math. Res. Not. 2015, No. 3, 651-687 (2015).




\bibitem{EvansBook2010}
Evans, L. C,
Partial differential equations. 2nd ed.
Graduate Studies in Mathematics 19. Providence, RI: American Mathematical Society. xxi, 749 p. (2010).



\bibitem{GilbargTrudinger2011}
Gilbarg, D, Trudinger, N. S,
Elliptic partial differential equations of second order. Reprint of the 1998 ed.
Classics in Mathematics. Berlin: Springer. xiii, 517 p. (2001).


\bibitem{GuoPeng2024}
Guo, Y, Peng, S,
Asymptotic behavior and classification of solutions to Hartree type equations with exponential nonlinearity.
J. Geom. Anal. 34, No. 1, Paper No. 23, 21 p. (2024).



\bibitem{HanLi2010}
Han, Z.-C, Li, Y,
On the local solvability of the Nirenberg problem on $\mathbb{S}^2$.
Discrete Contin. Dyn. Syst. 28, No. 2, 607-615 (2010).



\bibitem{HardyLittlewood1928}
Hardy, G. H, Littlewood, J. E,
Some properties of fractional integrals. I.
M. Z. 27, 565-606 (1928).




\bibitem{HebeyWei2013}
Hebey, E, Wei, J,
Schr\"odinger-Poisson systems in the 3-sphere.
Calc. Var. Partial Differ. Equ. 47, No. 1-2, 25-54 (2013).



\bibitem{Heinz1995}
Heinz, E,
\"Uber die Eindeutigkeit beim Cauchyschen Anfangswertproblem einer elliptischen Differentialgleichung zweiter Ordnung.
Nachr. Akad. Wiss. G?ttingen, Math.-Phys. Kl., Math.-Phys.-Chem. Abt. 1955, 1-12 (1955).


\bibitem{IanniRuiz2012}
Ianni, I, Ruiz, D,
Ground and bound states for a static Schr\"odinger-Poisson-Slater problem.
Commun. Contemp. Math. 14, No. 1, 1250003, 22 p. (2012).



\bibitem{LeiLiMa2012}
Lei, Y, Li, C, Ma, C,
Asymptotic radial symmetry and growth estimates of positive solutions to weighted Hardy-Littlewood-Sobolev system of integral equations.
Calc. Var. Partial Differ. Equ. 45, No. 1-2, 43-61 (2012).





\bibitem{Li1993}
Li, Y. Y,
On $-\Delta u=K(x)u^5$ in $\mathbb{R}^3$.
Commun. Pure Appl. Math. 46, No. 3, 303-340 (1993).



\bibitem{LiShafrir1994}
Li, Y, Shafrir, I,
Blow-up analysis for solutions of $-\Delta u=Ve^u$ in dimension two.
Indiana Univ. Math. J. 43, No. 4, 1255-1270 (1994).




\bibitem{Lieb1983}
Lieb, E. H,
Sharp constants in the Hardy-Littlewood-Sobolev and related inequalities.
Ann. Math. (2) 118, 349-374 (1983).



\bibitem{LiuMaXia2021}
Liu, X, Ma, S, Xia, J,
Multiple bound states of higher topological type for semi-classical Choquard equations.
Proc. R. Soc. Edinb., Sect. A, Math. 151, No. 1, 329-355 (2021).


\bibitem{MaZhao2010}
Ma, L, Zhao, L,
Classification of positive solitary solutions of the nonlinear Choquard equation.
Arch. Ration. Mech. Anal. 195, No. 2, 455-467 (2010).



\bibitem{MolicaRadulaescuSeradeiBook2016}
Molica, B. G, Radulescu, V. D, Servadei, R,
Variational methods for nonlocal fractional problems.
Encyclopedia of Mathematics and its Applications 162. Cambridge: Cambridge University Press. xvi, 383 p. (2016).



\bibitem{MorozVanSch2013}
Moroz, V, Van Schaftingen, J,
Groundstates of nonlinear Choquard equations: existence, qualitative properties and decay asymptotics.
J. Funct. Anal. 265, No. 2, 153-184 (2013).



\bibitem{MorozVanSch2015}
Moroz, V, Van Schaftingen, J,
Existence of groundstates for a class of nonlinear Choquard equations.
Trans. Am. Math. Soc. 367, No. 9, 6557-6579 (2015).



\bibitem{MorozVanS2017}
Moroz, V, Van Schaftingen, J,
A guide to the Choquard equation.
J. Fixed Point Theory Appl. 19, No. 1, 773-813 (2017).




\bibitem{Penrose1996}
Penrose, R,
On gravity's role in quantum state reduction.
Gen. Relativ. Gravitation 28, No. 5, 581-600 (1996).










\bibitem{PereraAgOReBook2010}
Perera, K, Agarwal, R. P, O'Regan, D,
Morse theoretic aspects of $p$-Laplacian type operators.
Mathematical Surveys and Monographs 161. Providence, RI: American Mathematical Society. xx, 141 p. (2010).



\bibitem{RabinowitzBook1986}
Rabinowitz, P. H,
Minimax methods in critical point theory with applications to differential equations.
Regional Conference Series in Mathematics 65. Providence, RI: American Mathematical Society. viii, 100 p. (1986).






\bibitem{Vaira2011}
Vaira, G,
Ground states for Schr\"odinger-Poisson type systems.
Ric. Mat. 60, No. 2, 263-297 (2011).


\bibitem{WangXia2020}
Wang, Z.-Q, Xia, J,
Saddle solutions for the Choquard equation. II.
Nonlinear Anal., Theory Methods Appl., Ser. A, Theory Methods 201, Article ID 112053, 25 p. (2020).



\bibitem{Watson1995}
Watson, G. N,
A treatise on the theory of Bessel functions. 2nd ed.
Cambridge Mathematical Library. Cambridge: Cambridge Univ. Press. vi, 804 p. (1995).



\bibitem{WeiWinter2009}
Wei, J, Winter, M,
Strongly interacting bumps for the Schr\"odinger-Newton equations.
J. Math. Phys. 50, No. 1, 012905, 22 p. (2009).



\bibitem{WillemBook1996}
Willem, M,
Minimax theorems.
Progress in Nonlinear Differential Equations and their Applications. 24. Boston: Birkh\"auser. viii, 159 p. (1996).




\bibitem{XiaWang2019}
Xia, J, Wang, Z.-Q,
Saddle solutions for the Choquard equation.
Calc. Var. Partial Differ. Equ. 58, No. 3, Paper No. 85, 30 p. (2019).




\end {thebibliography}
}
\end{document}